\documentclass{article}%
\usepackage{amsfonts}%
\usepackage{amsmath}%
\usepackage{amssymb}%
\usepackage{graphicx}
\begin{document}

\begin{center}
\textbf{JUCYS-MURPHY ELEMENTS}

\textbf{AND A SYMMETRIC FUNCTION IDENTITY}

\bigskip\ 

\textsc{Jennifer R. Galovich}

St. John's University

Collegeville, MN 56321

\medskip\ 
\end{center}

\bigskip\textsc{Abstract: }\ {\small Consider the elements of the group
algebra \ }$CS_{n}${\small \ \ given by }$R_{j}=\sum_{i=1}^{j-1}%
(ij),${\small \ for }$2\leq j\leq n.${\small \ Jucys [\textbf{3 - 5}] and
Murphy[\textbf{7}] showed that these elements act diagonally on elements of
}$S_{n}${\small \ and gave explicit formulas for the diagonal entries. \ We
give a new, combinatorial proof of this work in case }$j=n${\small \ and
present several similar results which arise from these combinatorial methods.}

\smallskip\ 

In a series of papers published early in the twentieth century, Alfred Young
described three forms for the irreducible representations of the symmetric
group $S_{n}$ [\textbf{12}]. Among these, the seminormal form enjoys several
nice properties:

(i) Matrices corresponding to adjacent transpositions can be computed
explicitly and easily.

(ii) The representation restricts from $S_{n}$ to $S_{n-1}$ in block diagonal
form with no change of basis required.

A.A. Jucys [\textbf{3 - 5}] and G.E. Murphy [\textbf{7}] gave a different
construction by introducing elements of the group algebra $\mathbb{C}S_{n}$
which act diagonally, and from which Young's seminormal form can be recovered.
Moreover, the diagonal entries of these Jucys-Murphy elements are easy to describe..

In 1994, the late S. Kerov asked for a combinatorial proof of a certain
symmetric function identity. [\textbf{1}]. \ As it happens, that identity is
equivalent to the action of a particular Jucys-Murphy element. \ In this note,
therefore, we first provide the requested combinatorial proof. \ Then we
present several variations on the Jucys-Murphy theme which are suggested by
these methods.

\smallskip\ 

We begin with some definitions and notation. A \textit{partition} of a
positive integer $n$ is a sequence $\lambda=(\lambda_{1},\lambda
_{2},...,\lambda_{k})$ where $\lambda_{1}\geq\lambda_{2}\geq\cdots\geq
\lambda_{k}>0$ and $\sum\lambda_{i}=n.$ We write $\lambda\vdash n.$ To each
$\lambda\vdash n$ we associate its \textit{Ferrers diagram} $F_{\lambda}$, a
left- and top-justified array of squares or cells. More precisely, if
$\lambda=(\lambda_{1},\lambda_{2},...,\lambda_{k})$ then $F_{\lambda}$
consists of $k$ left-justified rows of lengths $\lambda_{1},\lambda
_{2},...,\lambda_{k}$, reading from top to bottom. \ If $x$ is a cell in
$F_{\lambda}$ located in row $i$ and column $j$ we define $c\left(  x\right)
$, the \textit{content} of $x,$ by $c\left(  x\right)  =j-i.$(See Figure 1a.)

Note that by deleting any corner cell in $F_{\lambda}$ we automatically obtain
the Ferrers diagram of a partition of $n-1.$ Partitions obtained in this way
will be of particular interest; such a partition will be denoted $\lambda-x$,
where $x$ names the deleted corner cell. Similarly, by adding a cell $x$ to an
unoccupied corner of $F_{\lambda}$ we obtain a partition of $n+1$. A partition
which arises in this way will be denoted $\lambda+x$.

If $\lambda$ and $\mu$ are partitions we define the skew diagram $\lambda-\mu$
to be the set theoretic difference $F_{\lambda}\setminus$ $F_{\mu}$. A
\textit{rimhook} is a skew diagram which contains no $2\times2$ square as a
subset. The length of a rimhook is the number of cells it contains; the
height, $ht$, is one less than the number of rows it occupies.(See Figure 1b.)

\begin{center}
\bigskip$%
\begin{array}
[c]{cccc}%
\bullet & \bullet & \bullet & \bullet\\
\bullet & \bullet & \bullet & \\
\bullet & \bullet &  &
\end{array}
$

\smallskip\ 

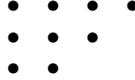
Figure 1a: $F_{\lambda}$ for $\lambda=(4,3,2)$

\bigskip

$\;%
\begin{array}
[c]{ccc}
& \bullet & \bullet\\
\bullet & \bullet & \\
\bullet &  &
\end{array}
$

\smallskip\ 

Figure 1b:

The rimhook $\lambda-\mu$ for $\lambda=(4,3,2)$ and $\mu=(2,1,1).$

Length of $\lambda-\mu$ is 5; $ht(\lambda-\mu)=2.$
\end{center}

Let $S_{n}$ denote the symmetric group on $n$ letters. Throughout, we will
write elements of $S_{n}$ using cycle notation, identifying the cycle type of
a permutation as the partition determined by the cycle lengths.

To each $\lambda\vdash n$ we also associate $\chi^{\lambda}$, the irreducible
character of $S_{n}$ corresponding to $\lambda$. As functions on $S_{n}$ the
$\chi^{\lambda}$s are constant on the conjugacy classes of $S_{n};$ indeed,
the collection $\{\chi^{\lambda}:\lambda\vdash n\}$ forms an orthonormal basis
for $CF_{n}$, the space of all class functions on $S_{n}.$

We use $\Lambda_{n}$ to denote the space of homogeneous symmetric functions of
degree $n$. Of the six standard bases for $\Lambda_{n}$ two are important for
the present work. They are $\{p_{\lambda}:\lambda\vdash n\}$ (the power sum
symmetric functions) and $\{s_{\lambda}:\lambda\vdash n\}$ (the Schur
functions), an orthonormal basis for $\Lambda_{n}.$ (We refer the reader to
[\textbf{6, 9, 10}] for more detailed information about symmetric functions.)

Our main theoretical tool is the (Frobenius) characteristic map which relates
the spaces $CF_{n}$ and $\Lambda_{n}$:

\smallskip\ 

\textsc{Definition:} Let $f\in CF_{n}.$ The characteristic map ch$^{n}%
:CF_{n}\rightarrow\Lambda_{n}$ is defined by%

\[
\text{ch}^{n}(f)=\sum_{\mu\vdash n}f(\mu)\cdot\frac{p_{\mu}}{z_{\mu}}
\]

\noindent where $p_{\mu}$ is the power sum symmetric function corresponding to
$\mu$ and $z_{\mu}=m_{1}!1^{m_{1}}m_{2}!2^{m_{2}}\cdots m_{k}!k^{m_{k}}.$

The map ch$^{n}\,\,$is an isometry; moreover, ch$^{n}(\chi^{\lambda
})=s_{\lambda}.$ [\textbf{9, }p. 163].

\smallskip\ 

Given any $\sigma\in S_{n}$ and $2\leq j\leq n$ we define the Jucys- Murphy
element $R_{j}(\sigma)\in\mathbb{C}S_{n}$ by%

\[
R_{j}(\sigma)=\sum_{i=1}^{j-1}\sigma\cdot(i\,j).
\]

\noindent The Jucys-Murphy elements $R_{j}$ have many interesting properties
and applications [\textbf{2}]. Chief among these is the fact that for any
$\lambda\vdash n$ the character%

\[
\chi^{\lambda}(R_{j}(\sigma))\equiv\sum_{i=1}^{j-1}\chi^{\lambda}(\sigma
\cdot(i\,j))
\]
can be computed easily in terms of the contents of the cells of $F_{\lambda} $
[\textbf{7}]. In case $j=n$, that theorem can be formulated as follows:

\smallskip\ 

\textbf{Theorem 1[Murphy, 7]. } \textit{Let }$\sigma\in S_{n}$\textit{\ with
}$\sigma(n)=n.$\textit{\ Let }$\overline{\sigma}$\textit{\ denote the
restriction of }$\sigma$\textit{\ to }$\{1,2,...n-1\}.$\textit{\ For any
}$\lambda\vdash n$\textit{\ we have}%

\begin{equation}
\chi^{\lambda}(R_{n}(\sigma))=\sum_{i=1}^{n-1}\chi^{\lambda}(\sigma
\cdot(i\,n))=\sum_{x}c(x)\cdot\chi^{\lambda-x}(\overline{\sigma})\qquad
\end{equation}

\textit{\noindent}\noindent\noindent\textit{where the sum on the right is
taken over all corner cells }$x$\textit{\ of }$F_{\lambda}$\textit{.}

Remark: A slightly different formulation of $(1)$ follows immediately from the
Branching Rule [\textbf{9}, p. 77]:%

\begin{equation}
\chi^{\lambda}(R_{n}(\sigma)+\chi^{\lambda}(\sigma)=\sum_{i=1}^{n}%
\chi^{\lambda}(\sigma\cdot(i\,n))=\sum_{x}(1+c(x))\cdot\chi^{\lambda
-x}(\overline{\sigma})\qquad
\end{equation}

\textit{Proof of Theorem 1. }Note that $\chi^{\lambda}(R_{n}(\sigma)$ depends
only on the cycle type $\mu$ of $\overline{\sigma}$, for $(\sigma\cdot(i\,n)$
has the same cycle structure as $\overline{\sigma}$ except that the cycle in
$\sigma$ containing $i$ is augmented by inserting $n$. \ For example, let
$\sigma=(253)(1)(4)\in S_{6}.$Then $\overline{\sigma}$= (253)(1)(4) and has
cycle type $\widehat{\mu}=1^{2}3^{1}.$ In this case
\begin{align*}
& R_{6=(2653)(1)(4)}+(2563)(1)(4)+(2536(1)(4)\\
& +(253)(16)(4)+(253)(1)(46).
\end{align*}
Of these five summands, two have type 1$^{1}2^{1}3^{1}$ and three have type
1$^{2}3^{0}4^{1}.$ In general, if the cycle type of $\overline{\sigma}$is
$\mu=1^{m_{1}}2^{m_{2}}\cdots k^{m_{k}}$ then among the $n-1=\sum j\cdot
m_{j}$ summands of $R_{n}(\sigma)$ there are $j\cdot m_{j}$ with cycle type
$\widehat{\mu}=1^{m_{1}}2^{m_{2}}\cdots j^{m_{j}-1}(j+1)^{m_{j}+1}\cdots
k^{m_{k}}.$

When $\widehat{\mu}$ is obtained from $\mu$ in this way, by replacing an
existing part of size $j$ with a part of size $j+1$, we write $\widehat{\mu
}>_{j}\mu.$ Thus%

\[
\sum\limits_{i=1}^{n-1}\chi^{\lambda}(\sigma\cdot(i\,n))=\sum\limits_{j\geq
1}\sum_{\widehat{\mu}>_{j}\mu}\chi^{\lambda}(\widehat{\mu})\cdot(j\cdot
m_{j}).
\]

\noindent Since the right side of $(1)$ can also be viewed as a function of
$\mu$, the key is to apply the Frobenius characteristic function to both sides.

Working first with the left side of $(1)$, let $F(\mu)=\sum\limits_{j\geq
1}\sum\limits_{\widehat{\mu}>_{j}\mu}\chi^{\lambda}(\widehat{\mu})\cdot(j\cdot
m_{j}).$ Then
\[
\begin{array}
[c]{lll}%
\text{ch}^{n}(F) & = & \sum_{\mu\vdash n-1}F(\mu)\cdot\frac{p_{\mu}}{z_{\mu}%
}\\
&  & \\
& = & \sum_{\mu\vdash n-1}\left[  \sum\limits_{j\geq1}\sum\limits_{\widehat
{\mu}>_{j}\mu}\chi^{\lambda}(\widehat{\mu})\cdot(j\cdot m_{j})\right]
\cdot\frac{p_{\mu}}{z_{\mu}}%
\end{array}
\]

\noindent Note that $\widehat{\mu}>_{j}\mu$ implies%

\[
\frac{j\cdot m_{j}}{z_{\mu}}=\frac{(j+1)(1+m_{j+1})}{z_{\widehat{\mu}}}
\]

\noindent and similarly, $p_{\mu}=p_{\widehat{\mu}}\cdot%
%TCIMACRO{\QDOVERD{.}{.}{p_{j}}{p_{j+1}}}%
%BeginExpansion
\genfrac{.}{.}{}{0}{p_{j}}{p_{j+1}}%
%EndExpansion
.$ Substituting and reversing the order of summation, we have%

\[
\begin{array}
[c]{cll}%
\text{ch}^{n}(F) & = & \sum_{\mu\vdash n-1}\sum\limits_{j\geq1}\sum\limits_{
\widehat{\mu}>_{j}\mu}\chi^{\lambda}(\widehat{\mu})\cdot\frac{(j+1)(1+m_{j+1}%
)}{z_{\widehat{\mu}}}\cdot p_{\widehat{\mu}}\cdot%
%TCIMACRO{\QDOVERD{.}{.}{p_{j}}{p_{j+1}}}%
%BeginExpansion
\genfrac{.}{.}{}{0}{p_{j}}{p_{j+1}}%
%EndExpansion
\\
&  & \\
& = & \sum\limits_{ \widehat{\mu}\vdash n}\chi^{\lambda}(\widehat{\mu}%
)\cdot\frac1{z_{\widehat{\mu}}}\cdot\left[  p_{\widehat{\mu}}\sum
\limits_{j\geq1}p_{j}(j+1)\frac{(1+m_{j+1})}{p_{j+1}}\right] \\
&  & \\
& = & \sum\limits_{\widehat{\mu}\vdash n}\chi^{\lambda}(\widehat{\mu}%
)\cdot\frac1{z_{\widehat{\mu}}}\cdot\left[  \sum\limits_{j\geq1}%
p_{j}(j+1)\frac\partial{\partial p_{j+1}}(p_{\widehat{\mu}})\right]
\end{array}
\]

\noindent since $1+m_{j+1}$ is the multiplicity of $j+1$ in $\widehat{\mu}.$
Continuing,
\[
\begin{array}
[c]{cll}%
\text{ch}^{n}(F) & = & \sum\limits_{j\geq1}p_{j}(j+1)\frac\partial{\partial
p_{j+1}}\left(  \sum\limits_{ \widehat{\mu}\vdash n}\chi^{\lambda}%
(\widehat{\mu})\cdot\frac{p_{\widehat{\mu}}}{z_{\widehat{\mu}}}\right) \\
&  & \\
& = & \sum\limits_{j\geq1}p_{j}(j+1)\frac\partial{\partial p_{j+1}}%
(s_{\lambda}).
\end{array}
\]

Returning now to the right side of $(1)$, we have%

\[
\text{ch}^{n}\text{ }\sum\limits_{x}\chi^{\lambda-x}\cdot c(x)=\sum
\limits_{x}s_{\lambda-x}\cdot c(x);
\]

\noindent thus it suffices to show that
\begin{equation}
\sum\limits_{j\geq1}p_{j}(j+1)\frac{\partial}{\partial p_{j+1}}(s_{\lambda
})=\sum\limits_{x}s_{\lambda-x}\cdot c(x).\qquad\qquad\qquad
\end{equation}

As operators on Schur functions, both $p_{j}$ and $Dp_{j}=j\frac
\partial{\partial p_{j}}$ can be interpreted in terms of rimhooks:

\smallskip\ 

\textbf{Lemma.} \textit{Let }$\lambda$ \textit{be a partition, }$j$\textit{\ a
positive integer. Then}

\smallskip\ 

(i) $p_{j}s_{\lambda}=\sum\limits_{\nu}(-1)^{ht(\nu-\lambda)}s_{\nu}$

(ii) $(Dp_{j})\cdot s_{\lambda}=\sum\limits_{\nu}(-1)^{ht(\lambda-\nu)}s_{\nu
}$

\textit{\noindent where the sums are taken over all partitions }$\nu
$\textit{\ such that }$\nu-\lambda$ \textit{(resp. }$\lambda-\nu)$\textit{\ is
a rimhook of length }$j$.

\textit{\smallskip\ }

\textit{Proof}. (i) [\textbf{6,} p. 31]

(ii) Since $\{s_{\lambda}\}$ is an orthonormal basis for $\Lambda_{n}$ it is
enough to compute $\langle D(p_{j})s_{\lambda},s_{\nu}\rangle$. Using part (i)
and the fact that, for any symmetric function $f$, the operator $D(f)$ is the
adjoint of multiplication by $f$ [\textbf{6,} p. 43] we have%
\[
\begin{array}
[c]{lll}%
\langle D(p_{j})s_{\lambda},s_{\nu}\rangle & = & \langle s_{\lambda}%
,p_{j}s_{\nu}\rangle\\
& = & \langle s_{\lambda},\sum\limits_{\zeta}(-1)^{ht(\zeta-\nu)}s_{\zeta
}\rangle\\
& = & \sum\limits_{\zeta}(-1)^{ht(\zeta-\nu)}\langle s_{\lambda},s_{\zeta
}\rangle\\
& = & (-1)^{ht(\lambda-\nu)}%
\end{array}
\]

\noindent where the sums are taken over all partitions $\zeta$ such that
$\zeta - \nu$ is a rimhook of length $j$. 

We use the Lemma to recast the left side of (2) as follows:%

\[
\begin{array}
[c]{lll}%
\sum\limits_{j\geq1}p_{j}(j+1)\frac\partial{\partial p_{j+1}}(s_{\lambda}) &
= & \sum\limits_{j\geq1}p_{j}D(p_{j+1})(s_{\lambda})\\
&  & \\
& = & \sum\limits_{j\geq1}p_{j}\sum\limits_{\nu}(-1)^{ht(\lambda-\nu)}s_{\nu
}\\
&  & \\
& = & \sum\limits_{j\geq1}\sum\limits_{\nu}(-1)^{ht(\lambda-\nu)}p_{j}s_{\nu
}\\
&  & \\
& = & \sum\limits_{j\geq1}\sum\limits_{\nu}(-1)^{ht(\lambda-\nu)}%
\sum\limits_{\zeta}(-1)^{ht(\zeta-\nu)}s_{\zeta}%
\end{array}
\]

\noindent where the first sum is taken over partitions $\nu$ such that
$\lambda-\nu$ is a rimhook of length $j+1$ and the second is over partitions
$\zeta$ such that $\zeta-\nu\,$ is a rimhook of length $j$.

To complete the proof of Theorem 1, we need to establish the identity%

\begin{equation}
\sum\limits_{j\geq1}\sum\limits_{\nu}(-1)^{ht(\lambda-\nu)}\sum\limits_{\zeta
}(-1)^{ht(\zeta-\nu)}s_{\zeta}=\sum\limits_{x}s_{\lambda-x}\cdot c(x).\;
\end{equation}

Note that the coefficient of each $s_{\zeta}$ arises by considering all
possible ways in which one may obtain the shape $F_{\zeta}$ ($\zeta\vdash
n-1$) by removing a rimhook of length $j+1$ from $F_{\lambda}$ to obtain a
shape $F_{\nu}$ then adding to $F_{\nu}$ a rimhook of length $j$. For example,
if $\lambda=332$ and $\zeta=322$ then the coefficient of $s_{\zeta}$ arises
from the following cases:%
\[
\begin{array}
[c]{lllll}%
j & \nu & ht(\lambda-\nu) & ht(\zeta-\nu) & \text{net contribution}\\
1 & (2,2,2) & 1 & 0 & -1\\
2 & (3,1,1) & 1 & 1 & +1\\
3 & (3,1) & 1 & 1 & +1
\end{array}
\]

Therefore the coefficient of $s_{322}$ on the left side of (4) is +1; note
that $322$ = $332-x$ where $x=(2,3),$ so $c(x)=+1$ is also the coefficient of
$s_{322}\,$ on the right, as predicted. However, if $\zeta=43$ then the
coefficient of $s_{\zeta}$is zero, since the cases $j=1$ ($\nu=33$) and $j=4$
($\nu=21$) give signs +1 and -1, respectively. This result is consistent with
the right side, since $\zeta=43$ is not of the form $\lambda-x$ for any $x$.
In the same way, the general argument divides into two cases:

Case (i). $\zeta=\lambda-x$ for some $x$: Suppose that $x=(p,q).$ Then $\zeta$
arises from $\lambda$ by removing a rimhook which either begins or ends with
the cell $(p,q)$. There are $q-1$ possibilities beginning with $(p,q)$ and any
rimhook added on must have the same sign, since $(p,q)$ is a corner in
$F_{\lambda}.$ There are $p-1$ ways to remove a rimhook ending at $(p,q)$ and
any rimhook added on must have the opposite sign. Therefore the net
contribution is $(q-1)-(p-1)=q-p=c(x).$

Case (ii). $\zeta\neq\lambda-x$ for any $x$: I claim that in all such
instances $\zeta$ arises in exactly two ways of opposite sign. Note first that
in these cases, both of the skew shapes $\lambda-\zeta$ and $\zeta-\lambda$
must be (non-empty) rimhooks since they are contained in the set of deleted or
added cells. For example, if $\lambda=(4,3,2,2)$ and $\zeta=(6,3,1)$ then
$\lambda-\zeta=$ $%
\begin{array}
[c]{r}%
\bullet\\
\bullet\bullet
\end{array}
$and $\zeta-\lambda=\bullet\bullet.$ Moreover, there are exactly two ways in
which $\zeta$ arises from removing and then adding a rimhook: For considering
the cells which connect $\lambda-\zeta$ and $\zeta-\lambda$ , either all are
removed and then replaced, or none of them, as illustrated in Figure 2.

\begin{center}
\bigskip$%
\begin{array}
[c]{cccccc}%
\bullet & \bullet & \otimes & \otimes & \times & \times\\
\bullet & \otimes & \otimes &  &  & \\
\bullet & \bigcirc &  &  &  & \\
\bigcirc & \bigcirc &  &  &  &
\end{array}
$ \smallskip\ 

$d=2;a=1;r=2$

\smallskip\ 

or

\smallskip\ 

$%
\begin{array}
[c]{cccccc}%
\bullet & \bullet & \bullet & \bullet & \times & \times\\
\bullet & \bullet & \bullet &  &  & \\
\bullet & \bigcirc &  &  &  & \\
\bigcirc & \bigcirc &  &  &  &
\end{array}
$ \smallskip\ 

$d=2;a=1;r=0$

\smallskip\ 

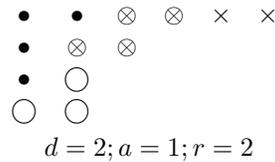
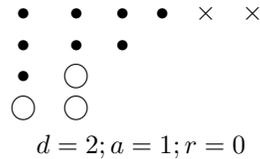
Figure 2: $\lambda=(4,3,2,2)$; $\zeta=(6,3,1)$

$\bullet$: original cell, unaffected

$\otimes$: original cell, removed and replaced

$\bigcirc$: original cell, removed but not replaced

$\times$: new cell
\end{center}

To show that these two ways have opposite signs, let $d(\lambda,\zeta)$ be the
number of rows in which some cells are deleted but not replaced; let
$a(\lambda,\zeta)$ be the number of rows in which some cells are added without
having been deleted; let $r(\lambda,\zeta)$ be the number of rows in which
cells are both deleted and replaced. If $r(\lambda,\zeta)=0$ then the sign
associated with $s_{\mu}$ is $(-1)^{a(\lambda,\zeta)+d(\lambda,\zeta)}.$
However, the connectedness of a rimhook guarantees that if $r(\lambda
,\zeta)\neq0$ then some row counted by $r(\lambda,\zeta)$ is also counted by
either $a(\lambda,\zeta)\,$ or $d(\lambda,\zeta).$ Therefore the sign in that
case is%
\[
(-1)^{a(\lambda,\zeta)+2r(\lambda,\zeta)+d(\lambda,\zeta)-1}=(-1)^{a(\lambda
,\zeta)+d(\lambda,\zeta)-1}.
\]

\noindent and the coefficient of $s_{\lambda}$ must be 0.

This completes the proof of Theorem 1. 

\bigskip

The same ideas used in the proof of Theorem 1 can also be used to produce
several interesting variations. In the first variation we replace the
Jucys-Murphy element $R_{n}$ by the analogous sum of 3-cycles. The summands on
the right side of equation (1) become values of characters corresponding to
partitions of $n-2$ obtained from $F_{\lambda}$ by removing two cells. Such a
partition will be denoted $\lambda-(x,y).$ The content $c(x,y)$ of a pair of
deleted cells is defined as follows:

(i) $c(x,y)=c(x)$ if $(x,y)$ forms a horizontal domino $
\begin{array}
[c]{cc}%
x & y
\end{array}
$

(ii) $c(x,y)=-c(x)$ if $(x,y)$ forms a vertical domino $
\begin{array}
[c]{c}%
x\\
y
\end{array}
$

(iii) $c(x,y)=-1$ if $x$ and $y$ are not contiguous.

\smallskip\ 

\textbf{Theorem 2:} \textit{Let }$\sigma\in S_{n}$\textit{\ with }%
$\sigma(n)=n$\textit{\ and }$\sigma(n-1)=n-1$\textit{, and define }%
$T_{n}(\sigma)=\sum_{i=1}^{n-2}\sigma\cdot(i\,n-1\,n).$\textit{\ Let
}$\overline{\sigma}$\textit{\ denote the restriction of }$\sigma$\textit{\ to
}$\{1,2,...n-2\}.$\textit{\ Then if }$\lambda$\textit{\ is any partition of
}$n$\textit{\ we have}%

\begin{equation}
\chi^{\lambda}(T_{n}(\sigma))\equiv\sum_{i=1}^{n-2}\chi^{\lambda}(\sigma
\cdot(i\,n-1\,n))=\sum_{(x,y)}\chi^{\lambda-(x,y)}(\overline{\sigma})\cdot
c(x,y)\;\;
\end{equation}

\smallskip\ 

\textit{Proof:} The Frobenius characteristic function can be applied to both
sides yielding%
\begin{equation}
\sum_{j\geq1}p_{j}Dp_{j+2}(s_{\lambda})=\sum_{(x,y)}s_{\lambda-(x,y)}\cdot
c(x,y)\quad
\end{equation}

The rest of the argument is essentially identical to that of Theorem 1 except
for the case in which $x$ and $y$ are not contiguous. That situation arises
exactly when $x$ and $y$ are the head and tail of a rimhook $\upsilon$ of
length $j+2$ replaced by the length $j$ rimhook $\nu-(x,y)$.\thinspace The
latter occupies one fewer row than $\nu$ so contributes the resulting shape
with multiplicity $-1.$

\smallskip\ 

A different variation on Theorem 1 is obtained by reversing the roles of $j$
and $j+1$ in equation (3). We have the following theorem; again the proof is
essentially the same as that of Theorem 1.

\smallskip\ 

\textbf{Theorem 3.} \textit{Let }$\lambda\vdash n-1$\textit{. Let }$\sigma\in
S_{n}$\textit{\ and set }$V_{n}(\sigma)=\sum\limits_{i\neq\sigma(i)}%
\sigma\cdot(i\;\sigma(i)).$\textit{\ When }$i\neq\sigma(i)$\textit{\ the
permutation }$\sigma\cdot(i\;\sigma(i))$\textit{\ has a fixed point and so may
be considered as an element of }$S_{n-1}$\textit{. With this in mind, we have}%

\begin{equation}
\chi^{\lambda}(V_{n}(\sigma))\equiv\sum_{i\neq\sigma(i)}\chi^{\lambda}%
(\sigma\cdot(i\;\sigma(i))=\sum_{x}\chi^{\lambda+x}(\sigma)\cdot c(x).\quad
\end{equation}

\smallskip

As in the case of Theorem 1, the Branching Rule immediately gives:

\
\begin{equation}
\sum_{i=1}^{n}\chi^{\lambda}(\sigma\cdot(i\;\sigma(i))=\sum_{x}\chi
^{\lambda+x}(\sigma)\cdot(1+c(x)).\qquad
\end{equation}

Remark: \ The Jucys- Murphy element $R_{n}$ acts on a permutation $\sigma$
which has a fixed point and, as noted earlier, the summands of $R_{n}(\sigma)$
are obtained by removing the fixed point and inserting it into each of the
cycles of $\overline{\sigma}$ in all possible ways. On the other hand, for an
arbitrary $\sigma\in S_{n}$, the element $V_{n}(\sigma)$ \textquotedblleft
inverses\textquotedblright\ this action by creating fixed points in all
possible ways.

Specializations of Theorem 3 lead to various corollaries. For example, if we
choose $\sigma$ to be the identity permutation, then the sum on the left in
equation (7) is empty and we have:%

\begin{equation}%
%TCIMACRO{\dsum \limits_{x}}%
%BeginExpansion
{\displaystyle\sum\limits_{x}}
%EndExpansion
f^{\lambda+x}\cdot c(x)=0
\end{equation}

where $f^{\lambda+x}$ is the dimension of the representation corresponding to
$\lambda+x.$

\smallskip\ 

Finally, the symmetric function approach we have described extends to the
hyperoctahedral group $B_{n}$ as well. \ One can use an analog of the
characteristic map described by John Stembridge [\textbf{11], }referring also
to the work of Arun Ram \textbf{8}] who has in fact extended the entire
Jucys-Murphy construction to types $B_{n}$, $D_{n}$, and $G_{2}.$

\begin{center}
\textsc{References}

\smallskip\ 
\end{center}

\bigskip1. S. Kerov, personal communication, 1994.

2. \ P. Diaconis and C. Greene, Applications of Murphy's elements,
unpublished, 1989.

\bigskip

3. \ A. A. Jucys, On the Young operators of symmetric groups,
\textit{Lithuanian Phys. J.} \textbf{6 }(1966) 163 - 180.\textbf{\ }

\bigskip

4. \ A. A. Jucys, Factorisation of Young's projective operators of symmetric
groups, \textit{Lithuanian Phys. J}. \textbf{11 }(1971) 1-10.

\bigskip

5. \ A. A. Jucys, Symmetric Polynomials and the center of the symmetric group
ring, \textit{Rep. Math. Phys} \textbf{5\ }(1974) 107-112.

\bigskip

6. I. G. MacDonald, Symmetric Functions and Hall Polynomials, Oxford
University Press, 1979.

\bigskip

7. \ G. E. Murphy, A new construction of Young's seminormal representation of
the symmetric groups, \textit{J. Algebra} \textbf{69} (1981) 287-297.

8. \ A Ram, Seminormal representations of Weyl groups and Iwahori-Hecke
algebras, \textit{Proc. London Math Soc. (3)} \textbf{75 }(1997) 99-133.

\smallskip\ 

9. \ B. Sagan, \textit{The Symmetric Group}, Wadsworth, 1991.

\smallskip\ \ 

10. \ R. Stanley, \textit{Enumerative Combinatorics}, Volume 2, Cambridge
University Press, 1999.

\bigskip

11. J. Stembridge, The projective representations of the hyperoctahedral
group, \textit{J. Algebra} \textbf{145} (1992) 396-453.

\bigskip

12. \ A. Young, Qualitative substitutional analysis III, \textit{Proc. London
Math. Soc (2)} \textbf{28} (1927) 255-292.

\smallskip\ 
\end{document}